\documentclass[12pt]{article}
\usepackage{amssymb}
\usepackage{amsmath}
\usepackage{array}
\usepackage[russian]{babel}
\usepackage[pdftex]{graphics}

\usepackage[X2,T2A]{fontenc}
\usepackage[cp1251]{inputenc}

\DeclareMathOperator{\RRe}{Re} 
\DeclareMathOperator{\mmod}{mod}

\DeclareMathOperator{\vep}{\varepsilon}

\newcommand{\prsum}{\mathop{{\sum}'}}

\multlinegap=0em \DeclareFontFamily{T1}{msb}{}
\DeclareFontShape{T1}{msb}{m}{ol}{<5> <6> <7> <8> <9> gen * msbm
<10> <10.95> <12> <14.4> <17.28> <20.74> <24.88> msbm10}{}
\DeclareSymbolFont{AMSb}{T1}{msb}{m}{ol} \multlinegap=0em

\renewcommand{\S}{\mathhexbox278}
\renewcommand{\le}{\operatorname{\leqslant}}
\renewcommand{\ge}{\operatorname{\geqslant}}

\begin{document}

\begin{flushleft}
MSC 11L05
\end{flushleft}

\begin{center}
{\rmfamily\bfseries\normalsize On Kloosterman sums with multiplicative coefficients}
\end{center}

\begin{center}
{\rmfamily\bfseries\normalsize M.A.~Korolev\footnote{This work is supported by the Russian Science Foundation under grant 14-11-00433
and performed in Steklov Mathematical Institure of Russian Academy of Sciences.}}
\end{center}

\vspace{0.5cm}

\fontsize{11}{12pt}\selectfont

\textbf{Abstract.} The series of some new estimates for the sums of the type
\[
S_{q}(x;f)\,=\,\prsum\limits_{n\le x}f(n)e_{q}(an^{*}+bn)
\]
is obtained. Here $q$ is a sufficiently large integer, $\sqrt{q}(\log{q})\!\ll\!x\le q$,
$a,b$ are integers, $(a,q)=1$, $e_{q}(v) = e^{2\pi iv/q}$, $f(n)$
is a multiplicative function, $nn^{*}\equiv 1 \pmod{q}$ and the prime sign means that $(n,q)=1$.
These estimates improve the previous results of such type belonging to K.~Gong and C.~Jia.

\vspace{0.2cm}


\fontsize{12}{15pt}\selectfont

\vspace{0.5cm}

\textbf{\S 1. Introduction}

Here we prove a series of some new estimates for incomplete weighted Kloosterman sums, that is,
for the sums of the following type:
\begin{equation}\label{lab_01}
S_{q}(x;f)\,=\,\prsum\limits_{n\le x}f(n)e_{q}(an^{*}+bn).
\end{equation}
Here $q$ is a sufficiently large integer, $\sqrt{q}\,(\log{q})\!\ll\!x\le q$, is any fixed number, $a,b$ are integers, $(a,q)=1$, $e_{q}(v) = e^{2\pi iv/q}$, $f(n)$
is a multiplicative function and the prime sign means that $(n,q)=1$. The symbol $n^{*}$ stands for the inverse residue for $n \pmod q$, that is, for the solution of the congruence $nn^{*}\equiv 1 \pmod{q}$.

The sums (\ref{lab_01}) with M\"{o}bius function $f(n) = \mu(n)$ were studied in \cite{Deng_1999}, \cite{Hajela_Pollington_Smith_1988}, \cite{Wang_Zheng_1998} and the following estimate was obtained in two last papers:
\begin{equation}\label{lab_02}
S_{q}(x;f)\,\ll\,x\tau(q)\bigl(q^{-1/2}(\log{x})^{5/2}\,+\,q^{1/5}x^{-1/5}(\log{x})^{13/5}\bigr),
\end{equation}
where, as usual, $\tau_{k}(q)$ is the divisor function, $\tau(q) = \tau_{2}(q)$. Since $\tau(q)\ll q^{\vep}$ for any fixed $\vep>0$, one can check that
the estimate (\ref{lab_02}) is non-trivial if
\[
q^{1+\vep}\ll x\ll \exp{\bigl(q^{1/5-\vep}\bigr)}.
\]
Recently, K.~Gong and C.~Jia \cite{Gong_Jia_2016} proved that in the case $b\equiv 0\pmod {q}$ the inequality
\begin{equation}\label{lab_03}
S_{q}(x;f)\,\ll\,x\bigl(\bigl(\tau(q)q^{-1}\log\log{q}\bigr)^{1/2}\log\log{x}\,+
\,q^{1/4+\vep}x^{-1/2}(\log{x})^{1/2}\,+\,(\log\log{x})^{-1/2}\bigr)
\end{equation}
holds for any multiplicative function $f(n)$ satisfying the condition $|f(n)|\le 1$. Obviously, this bound is non-trivial when
\[
q^{1/2+\delta_{1}}\ll x \ll \exp{\bigl(\exp{\bigl(q^{1/2-\delta_{2}}\bigr)}\bigr)}
\]
for some positive $\delta_{j} = \delta_{j}(\vep)>0$.

The aim of this paper is to improve a general bound (\ref{lab_03}) and to obtain some particular estimates for $S_{q}(x;f)$
for $\sqrt{q}\,(\log{q})\!\ll\!x\le q$ and for some particular functions $f(n)$. The main results of the paper are the following
(in theorems 1-4, $f(n)$ denotes any multiplicative function satisfying the condition $|f(n)|\le 1$).

\vspace{0.3cm}

\textsc{Theorem 1.} \emph{Let $0<\vep < 0.5$ be any fixed constant, $q\ge q_{1}(\vep)$ is a sufficiently large integer, and suppose that
$q^{1/2+\vep}\ll x\le q$. Then the following estimate holds:}
\[
|S_{q}(x;f)|\,\le\, 562\,x\,\frac{\log\log{q}}{\vep\log{q}}.
\]

\textsc{Theorem 2.} \emph{Let $\gamma>0$ be any fixed constant; $q\ge q_{2}(\gamma)$ is a sufficiently large integer, and suppose that}
\[
\sqrt{q}\,e^{(\log\log{q})^{1+\gamma}}\le x\le q.
\]
\emph{If $\tau(q)\le e^{0.25(\log\log{q})^{1+\gamma}}$ then the following estimate holds:}
\[
|S_{q}(x;f)|\,\le\, x\,\frac{562(2+\gamma)}{(\log\log{q})^{1+\gamma}}.
\]

\textsc{Theorem 3.} \emph{Let $\gamma>0$ be any fixed constant, $q\ge q_{3}(\gamma)$ is a sufficiently large integer, and suppose that}
\[
\sqrt{q}\,(\log{q})^{1+2\gamma}\le x\le q.
\]
\emph{If $\tau(q)\le (\log{q})^{\gamma/4}$ then the following estimate holds:}
\[
|S_{q}(x;f)|\,\le\, 281\,x\,\frac{\log\log\log{q}}{\gamma\log\log{q}}.
\]

\textsc{Theorem 4.} \emph{Let $\gamma>0$ be any fixed constant, $q\ge q_{4}(\gamma)$ is a sufficiently large prime, and suppose that}
\[
x\ge \sqrt{q}(\log{q})e^{(\log\log\log{q})^{1+\gamma}}.
\]
\emph{Then the following estimate holds:}
\[
|S_{q}(x;f)|\,\le \,\frac{562\,x}{(\log\log\log{q})^{\gamma\mathstrut}}.
\]

In some cases, the sum $S_{q}(x;f)$ is estimated with power-saving factor. Such estimates are based on
the bounds for double Kloosterman sums obtained by J.~Bourgain \cite[Appendix]{Bourgain_2005}.
In particular, the following assertions hold true.

\vspace{0.3cm}

\textsc{Theorem 5.} \emph{Let $0<\vep<0.1$ be any fixed constant, $q\ge q_{4}(\vep)$ is prime and suppose that
$q^{1/2+\vep}\le x\le q$. Then}
\[
\sum\limits_{n\le x}\mu(n)e_{q}(an^{*}+bn)\,\ll_{\vep}\, xq^{-c\vep^{4}}
\]
\emph{for some absolute constant} $c>0$.

\vspace{0.3cm}

\textsc{Theorem 6.} \emph{Let $0<\vep<0.1$ be any fixed constant, $k\ge 1$ is a fixed integer, $q\ge q_{5}(\vep;k)$ is prime and suppose that $q^{1/2+\vep}\le x\le q$. Then}
\[
\sum\limits_{n\le x}\tau_{k}(n)e_{q}(an^{*}+bn)\,\ll_{k,\vep}\, xq^{-c\vep^{4}}
\]
\emph{for some constant} $c = c(k)>0$.

\vspace{0.3cm}

\textsc{Acknowledgements.} This work was stimulated by fruitful discussions with pro\-fes\-sors Chaohua Jia
(Institute of Mathematics, Academia Sinica) and Ke Gong (Henan University) in Academy of Mathematics and System Science,
CAS (Beijing) during ``Chinese--Russian workshop of exponential sums and sumsets'' (November, 2015).
The author is warmly grateful to them for the support and hospitality.

\vspace{0.5cm}

\textbf{\S 2. Auxiliary assertions}

\vspace{0.5cm}

The below assertions are necessary for the proof of theorems 1--6.

\vspace{0.3cm}

\textsc{Lemma 1.} \emph{Suppose that $2\le y\le x$ and let $\Phi(x,y)$ be the quantity of numbers $n\le x$ free of prime divisors $\le y$. Then}
\[
\Phi(x,y)\,\le\,\frac{x}{\log{y}}\,+\,\frac{13.5x}{(\log{y})^{2}}.
\]

This estimate can be derived by standard technic of Selberg's sieve. The details of the proof are contained in \cite{Korolev_2016a}.

\vspace{0.3cm}

\textsc{Lemma 2.} \emph{Suppose that $15\le y<x$ and let $\Psi(x,y)$ be the quantity of numbers $n\le x$ free of prime divisors $>y$. Then}
\[
\Psi(x,y)\,\le\,Cxe^{-u/2},\quad u\,=\,\frac{\log{x}}{\log{y}},\quad C\,=\,67.21.
\]

\textsc{Proof.} We will follow the proof of Theorem 1 from \cite[Part III, ch. III.5, \S 5.1]{Tenenbaum_1995}. Setting $\alpha = \tfrac{2}{3}(\log{y})^{-1}$, we get
\[
\Psi(x,y)\,\le\,x^{3/4}+x^{-3\alpha/4}\sum\limits_{n\le x}f(n),
\]
where
\begin{equation*}
f(n)\,=\,\\
\begin{cases}
n^{\alpha}, & \text{if}\;\;P(n)\le y,\\
0, & \text{otherwise},
\end{cases}
\end{equation*}
and $P(n)$ denotes the largest prime divisor of $n$. Next, for any prime $p\le y$ and any $\nu\ge 1$ we obviously have
\[
f(p^{\nu})\,=\,p^{\alpha\nu}\,=\,\exp{\biggl(\frac{2\nu}{3}\,\frac{\log{p}}{\log{y}}\biggr)}\,\le\,e^{2\nu/3}.
\]
Hence, using Theorem 9 from \cite{Rosser_Schoenfeld_1962}, we have for any $z\ge 2$:
\[
\sum\limits_{p\le z}f(p)\log{p}\,\le\,e^{2/3}\sum\limits_{p\le \min{(y,z)}}\log{p}\,<\,A\min{(y,z)}\,\le\,Az,
\]
where $A = 1.01624e^{2/3}$. Next, the inequalities
\[
\alpha\,\le\,\frac{2}{3\log{15}}\,<\,\frac{1}{4}
\]
yield that
\[
\sum\limits_{p,\nu\ge 2}\frac{f(p^{\nu})}{p^{\nu\mathstrut}}\,\log{p^{\nu}}\,\le\,\sum\limits_{p\le y}\sum\limits_{\nu\ge 2}
\frac{\nu\log{p}}{p^{3\nu/4}}\,\le\,\sum\limits_{p}\frac{\log{p}}{p^{3/2}}\,\frac{2-p^{-3/4}}{(1-p^{-3/4})^{2\mathstrut}}\,<\,B,
\]
where $B=5.28475$. Thus, Theorem 5 from \cite[Part III, ch. III.3, \S 3.5]{Tenenbaum_1995} implies that
\[
\sum\limits_{n\le x}f(n)\,\le\,(A+B+1)\,\frac{x}{\log{x}}\sum\limits_{n\le x}\frac{f(n)}{n}.
\]
Since
\[
\frac{1}{\log{x}}\,\le\,\frac{1}{\log{y}}\,\le\,\prod\limits_{p\le y}\biggl(1-\frac{1}{p}\biggr),
\]
then
\begin{multline*}
\frac{1}{\log{x}}\sum\limits_{n\le x}\frac{f(n)}{n}\,\le\,\frac{1}{\log{x}}\prod\limits_{p\le y}\sum\limits_{\nu\ge 0}\frac{f(p^{\nu})}{p^{\nu\mathstrut}}\,\le\,\prod\limits_{p\le y}\biggl(1-\frac{1}{p}\biggr)\biggl(1-\frac{1}{p^{1-\alpha}}\biggr)^{-1}\,=\\
=\,\prod\limits_{p\le y}\biggl(1\,+\,\frac{p^{\alpha}-1}{p}\biggl(1-\frac{1}{p^{1-\alpha}}\biggr)^{-1}\biggr).
\end{multline*}
Now one can check that
\begin{multline*}
\biggl(1-\frac{1}{p^{1-\alpha}}\biggr)^{-1}\,\le\,\biggl(1-\frac{1}{2^{1-\alpha}}\biggr)^{-1}\,<\,\bigl(1-2^{-3/4}\bigr)^{-1}\,<\,\frac{5}{2},\\
p^{\alpha}-1\,<\,\alpha(\log{p})\bigl(1+0.5\alpha p^{\alpha}\bigr)\,\le\,\alpha(\log{p})\biggl(1\,+\,\frac{e^{2/3}}{3\log{y}}\biggr)\,<\,\frac{5}{4}\,\alpha\log{p}
\end{multline*}
for any prime $p\le y$. Using the corollary of Theorem 6 from \cite{Rosser_Schoenfeld_1962}, we find:
\begin{multline*}
\frac{1}{\log{x}}\sum\limits_{n\le x}\frac{f(n)}{n}\,\le\,\prod\limits_{p\le y}\biggl(1\,+\,\frac{25}{8}\,\alpha\,\frac{\log{p}}{p}\biggr)
\,<\,\exp{\biggl(\frac{25}{8}\,\alpha\sum\limits_{p\le y}\frac{\log{p}}{p}\biggr)}\,\le\\
\le\,\exp{\biggl(\frac{25}{8}\,\alpha\log{y}\biggr)}\,=\,e^{25/12}.
\end{multline*}
Thus we obtain:
\begin{multline*}
\sum\limits_{n\le x}f(n)\,\le\,Dx,\quad D = (A+B+1)e^{25/12},\\
\Psi(x,y)\,\le\,x^{3/4}+Dx^{1-3\alpha/4}\,=\,x^{1-3\alpha/4}\bigl(D\,+\,x^{-1/4+3\alpha/4}\bigr).
\end{multline*}
Since
\[
-\frac{1}{4}+\frac{3}{4}\,\alpha\,\le\,-\frac{1}{4}+\frac{1}{2\log{15}}<-\frac{1}{15.3},
\]
we finally get:
\[
\Psi(x,y)\,\le\,x^{1-3\alpha/4}\bigl(D\,+\,15^{-1/15.3}\bigr)\,<\,Cx\exp{\biggl(-\,\frac{1}{2}\,\frac{\log{x}}{\log{y}}\biggr)},
\]
where $C=67.21$. $\Box$

\vspace{0.3cm}

\textsc{Lemma 3.} \emph{Let $q\ge 3$, $a,b$, $M,N$ be the integers, $1<N<q$. Then the following estimates hold:}
\begin{multline*}
\biggl|\prsum\limits_{n=1}^{q}e_{q}(an^{*}+bn)\biggr|\,\le\,\tau(q)(a,b,q)^{1/2}q^{1/2},\\
\biggl|\prsum\limits_{M<n\le M+N}e_{q}(an^{*}+bn)\biggr|\,\le\,2\tau(q)(a,q)^{1/2}q^{1/2}\log{q}.
\end{multline*}

\vspace{0.3cm}

\textsc{Proof.} The derivation of the first estimate from the classical A.~Weil's theorem \cite{Weil_1948} is
contained in the paper of T.~Estermann \cite{Estermann_1961}. The second estimate easily follows from the chain of
relations
\begin{multline*}
\prsum\limits_{M<n\le M+N}e_{q}(an^{*}+bn)\,=\,\prsum\limits_{n=1}^{q}\biggl(\frac{1}{q}\sum\limits_{-q/2<c\le q/2}
\sum\limits_{M<m\le M+N}e_{q}(c(n-m))\biggr)e_{q}(an^{*}+bn)\,=\\
=\,\sum\limits_{-q/2<c\le q/2}\frac{\gamma_{c}}{|c|+1}\prsum\limits_{n=1}^{q}e_{q}(an^{*}+(b+c)n),\quad
\gamma_{c}=\frac{|c|+1}{q}\biggl(\sum\limits_{M<m\le M+N}e_{q}(-cm)\biggr)
\end{multline*}
and the obvious inequalities $|\gamma_{c}|\le 1$, $(a,b+c,q)\le (a,q)$ for any $c$, $-q/2<c\le q/2$.

\vspace{0.3cm}

\textsc{Lemma 4.} \emph{If $X\to +\infty$ then}
\[
\prod\limits_{p\le X}\biggl(1-\frac{1}{p}\biggr)^{-1}\,=\,e^{\gamma}(\log{X})\bigl(1\,+\,O\bigl(e^{-c\sqrt{\log{X}\mathstrut}}\bigr)\bigr),\quad c>0,
\]
\emph{where $\gamma$ denotes Euler's constant.}

\vspace{0.3cm}

This is the consequence from prime number theorem (see, for example, \cite[ch. 3, \S 5]{Prachar_1957}).

\vspace{0.5cm}

\textbf{\S 3. Proofs of theorems 1--4}

\vspace{0.5cm}

First we prove the main lemma which allows one to derive theorems 1-4 in the uniform manner.

\vspace{0.3cm}

\textsc{Lemma 5.} \emph{Suppose that}
\[
x\le q,\quad \frac{x}{\sqrt{q}(\log{q})(\log\log{q})^{2\mathstrut}}\,\to +\infty,
\]
\emph{and let $f(n)$ be any multiplicative function such that $|f(n)|\le 1$. Then the sum}
\[
S_{q}(x;f)\,=\,\prsum\limits_{n\le x}f(n)e_{q}(an^{*}+bn)
\]
\emph{satisfies the estimate $|S|\le x(\Delta_{1}+\Delta_{2}+\Delta_{3})$, where}
\[
\Delta_{1}\,=\,(C+3)\,\frac{\log{X}}{\log{Y}},\quad \Delta_{2}\,=\,\frac{7\log\log{q}}{\sqrt{X\log{X\mathstrut}}},\,\quad
\Delta_{3}\,=\,\frac{7}{2}\,\tau(q)\,\biggl(\frac{Y}{x}\,\sqrt{q}(\log{q})(\log\log{q})^{2}\biggr)^{\!1/2},
\]
\emph{$C$ is the constant from lemma 2, $X$ and $Y$ satisfy the following conditions}
\begin{equation}\label{lab_04}
15<X_{0}<X<0.5Y,\quad Y\le \frac{x}{\sqrt{q}(\log{q})(\log\log{q})^{2\mathstrut}}.
\end{equation}

\textsc{Proof.} Let us choose any $X, Y$ satisfying (\ref{lab_04}) and put $\mathcal{I} = (X,Y]$. Further, we denote denote by
$A_{r}$ the set of $n$, $1\le n\le x$, $(n,q)=1$, having exactly $r$ prime divisors
from $\mathcal{I}$ and counting with multiplicities ($r = 0,1,2,\ldots$).

Then any $n\in A_{0}$ is expressed in the form $n = uv$ or in the form $n = u$ where all prime divisors of $u$ does not exceed $X$
(or $u=1$) and all prime divisors of $v$ are greater than $Y$. Fixing $u$, we have at most $\Phi(xu^{-1},Y)$ possibilities for the factor $v$.
Summing over $u$ and using lemmas 1,4 we obtain that the number of $n = uv$, $n\in A_{0}$, does not exceed
\begin{multline*}
N_{1}\,=\,\sum\limits_{u}\Phi\biggl(\frac{x}{u},Y\biggr)\,\le\,\frac{2x}{\log{Y}}\sum\limits_{u}\frac{1}{u}\,\le\,
\frac{2x}{\log{Y}}\prod\limits_{p\le X}\biggl(1+\frac{1}{p}+\frac{1}{p^{2}}+\ldots\biggr)\,=\\
=\,\frac{2x}{\log{Y}}\prod\limits_{p\le X}\biggl(1-\frac{1}{p}\biggr)^{-1}<\,3x\,\frac{\log{X}}{\log{Y}}.
\end{multline*}
Next, using the notations of lemma 2, we conclude that the number of $n = u$, $n\in A_{0}$ is less than
\[
N_{2}\,=\,\Psi(x,X)\,\le \,Cx\exp{\biggl(-\,\frac{1}{2}\,\frac{\log{x}}{\log{X}}\biggr)},\quad C=67.21.
\]
In view of (\ref{lab_04}), we have $(\log{x})/(2\log{Y})\ge 1$. Thus we get:
\[
N_{2}\,\le\,Cx\exp{\biggl(-\,\frac{1}{2}\,\frac{\log{x}}{\log{Y}}\,\frac{\log{Y}}{\log{X}}\biggr)}\,\le\,Cx\exp{\biggl(-\,\frac{\log{Y}}{\log{X}}\biggr)}\,<\,Cx\,\frac{\log{X}}{\log{Y}}.
\]
Hence, setting
\[
S_{r}\,=\,\sum\limits_{n\in A_{r}}f(n)e_{q}(an^{*}+bn),
\]
we get
\[
|S_{0}|\,\le\,|A_{0}|\,\le\,N_{1}+N_{2}\,<\,(C+3)x\,\frac{\log{X}}{\log{Y}}.
\]
Suppose now that $r\ge 1$ and consider all the products $pm$ where $p$ and $m$ run independently the sets $\mathcal{I}$ and $A_{r-1}$. If $n\in A_{r}$ is not divisible
by squares of primes from $\mathcal{I}$ then it arises exactly $r$ times among these products with the conditions $(p,m)=1$, $(pm,q)=1$. Therefore,
\[
S_{r}\,=\,\frac{1}{r}\prsum\limits_{X<p\le Y}\sum\limits_{\substack{m\in A_{r-1} \\ pm\le x,\; (p,m)=1}}f(p)f(m)e_{q}(ap^{*}m^{*}+bpm)\,+\,\theta s_{r},
\]
where $|\theta|\le 1$ and $s_{r}$ denotes the number of $n\in A_{r}$ divisible by $p^{2}$ for some $p\in \mathcal{I}$. If we omit the condition $(p,m)=1$ in the inner sum
then the corresponding error is less than
\[
\frac{1}{r}\sum\limits_{p>X}\frac{x}{p^{2}}\,<\,\frac{xX^{-1}}{r}\,
\]
in absolute value. Thus we obtain:
\[
S_{r}\,=\,\frac{1}{r}\prsum\limits_{X<p\le Y}f(p)\sum\limits_{m\in A_{r-1}}f(m)e_{q}(ap^{*}m^{*}+bpm)\,+\,\theta_{1}\biggl(s_{r}+\frac{xX^{-1}}{r}\biggr).
\]
Now let us split $\mathcal{I}$ into intervals $Q<p\le Q_{1}$ where $Q_{1}\le 2Q$ and denote
\[
S_{r}(Q)\,=\,\sum\limits_{m\in A_{r-1}}f(m)\prsum\limits_{\substack{Q<p\le Q_{1} \\ pm\le x}}f(p)e_{q}(ap^{*}m^{*}+bpm).
\]
Since $m\le xQ^{-1}$ for any $m$ in the sum then
\[
|S_{r}(Q)|\,\le\,\prsum\limits_{m\le xQ^{-1}}\biggl|\prsum\limits_{\substack{Q<p\le Q_{1} \\ p\le x/m}}f(p)e_{q}(ap^{*}m^{*}+bpm)\biggr|.
\]
By Cauchy's inequality,
\begin{multline*}
|S_{r}(Q)|^{2}\,\le\,xQ^{-1}\prsum\limits_{m\le xQ^{-1}}\prsum\limits_{\substack{Q<p_{1}, p_{2}\le Q_{1} \\ p_{1},p_{2}\le x/m}}f(p_{1})\overline{f}(p_{2})
e_{q}(a(p_{1}^{*}-p_{2}^{*})m^{*}+b(p_{1}-p_{2})m)\,=\\
=\,xQ^{-1}\biggl(\;\prsum\limits_{m\le xQ^{-1}}\prsum\limits_{\substack{Q<p\le Q_{1} \\ p\le x/m}}|f(p)|^{2}\,+\,2\RRe
\prsum\limits_{Q<p_{1}<p_{2}\le Q_{1}}f(p_{1})\overline{f}(p_{2})\prsum\limits_{1\le m\le M}e_{q}(a_{1}m^{*}+b_{1}m)\biggr)\,\le\\
\le\,xQ^{-1}\biggl(xQ^{-1}\pi(Q_{1})\,+\,2\prsum\limits_{Q<p_{1}<p_{2}\le Q_{1}}\biggl|\prsum\limits_{1\le m\le M}e_{q}(a_{1}m^{*}+b_{1}m)\biggr|\biggr),
\end{multline*}
where $M = xp_{2}^{-1}$, $a_{1}=a(p_{1}^{*}-p_{2}^{*})$, $b_{1}=b(p_{1}-p_{2})$. Since $(p_{1}p_{2},q)=1$ then
\[
p_{1}p_{2}(p_{1}^{*}-p_{2}^{*})\equiv p_{2}-p_{1}\pmod{q}
\]
and the numbers $a_{1}, p_{2}-p_{1}$ has the same greatest common divisor with $q$: $(a_{1},q) = (p_{2}-p_{1},q) = \delta$. Using lemma 3, we find
\begin{multline*}
\biggl|\prsum\limits_{1\le m\le M}e_{1}(a_{1}m^{*}+b_{1}m)\biggr|\,\le\,2\tau(q)\sqrt{q}(\log{q})\sqrt{\delta},\\
\bigl|S_{r}(Q)\bigr|^{2}\,\le\,\frac{2x^{2}}{Q\log{Q}}\,+\,\frac{2x}{Q}\cdot 2\tau(q)\sqrt{q}(\log{q})
\prsum\limits_{\substack{Q<p_{1}<p_{2}\le Q_{1} \\ (p_{2}-p_{1},q)=\delta}}\sqrt{\delta}.
\end{multline*}
Obviously, $\delta\le p_{2}-p_{1}\le Q$. Hence, the last sum does not exceed
\begin{multline*}
\sum\limits_{\delta|q,\;\delta\le Q}\sqrt{\delta}\sum\limits_{Q<p_{1}\le Q_{1}}\sum\limits_{\substack{p_{2}\equiv p_{1}(\mmod\delta) \\ Q<p_{2}\le Q_{1}}}1\,\le\\
\le\,\sum\limits_{\delta|q,\;\delta\le Q}\sqrt{\delta}\sum\limits_{Q<p_{1}\le Q_{1}}\biggl(\frac{Q}{\delta}+1\biggr)\,\le\,2Q\sum\limits_{\delta|q,\;\delta\le Q}\frac{1}{\sqrt{\delta}}
\sum\limits_{Q<p_{1}\le Q_{1}}1\,<\\
<\,2Q\pi(Q_{1})\sum\limits_{\delta|q}\frac{1}{\sqrt{\delta}}\,<\,\frac{5Q^{2}\tau(q)}{\log{Q}}.
\end{multline*}
Therefore,
\[
\bigl|S_{r}(Q)\bigr|^{2}\,\le\,\frac{2x^{2}}{Q\log{Q}}\,+\,\frac{20xQ(\log{q})}{\log{Q}}\tau^{2}(q)\sqrt{q}\,<\,
\frac{2x^{2}}{Q\log{Q}}\,+\,xQ(\log{q})\tau^{2}(q)\sqrt{q}
\]
and hence
\begin{equation}\label{lab_05}
\bigl|S_{r}(Q)\bigr|\,\le\,\frac{1.5x}{\sqrt{Q\log{Q}\mathstrut}}\,+\,\tau(q)\bigl(xQ\sqrt{q}\log{q}\bigr)^{\!1/2}\,<
\,x\biggl(\frac{1.5}{\sqrt{Q\log{Q}\mathstrut}}\,+\,\tau(q)\biggl(\frac{Q\sqrt{q}\log{q}}{x}\biggr)^{\!\! 1/2}\,\biggr).
\end{equation}
Taking $Q = 2^{k}X$, $k = 0,1,2,\ldots, k_{0}$, where $2^{k_{0}}X\le Y< 2^{k_{0}+1}X$ and summing both parts of (\ref{lab_05})
over $k$, we obtain:
\begin{multline*}
\sum\limits_{k = 0}^{k_{0}}\bigl|S_{r}(Q)\bigr|\,<\,x\biggl(1.5\sum\limits_{k = 0}^{k_{0}}\frac{2^{-k/2}}{\sqrt{X\log{X}\mathstrut}}
\,+\,\tau(q)\biggl(\frac{\sqrt{q}\log{q}}{x}\biggr)^{\!\! 1/2}\sum\limits_{k = 0}^{k_{0}}(2^{k}X)^{1/2}\,\biggr)\,<\\
<\,\frac{7x}{2}\biggl(\frac{1.5}{\sqrt{X\log{X}\mathstrut}}\,+\,\tau(q)\biggl(\frac{Y\sqrt{q}\log{q}}{x}\biggr)^{\!\! 1/2}\,\biggr),\\
|S_{r}|\,<\,\frac{7x}{2r}\biggl(\frac{1.5}{\sqrt{X\log{X}\mathstrut}}\,+\,\tau(q)\biggl(\frac{Y\sqrt{q}\log{q}}{x}\biggr)^{\!\! 1/2}\,\biggr)+s_{r}+\frac{xX^{-1}}{r}.
\end{multline*}
Summing over $r$, $1\le r\le (\log{x})/\log{X}$ and noting that
\[
\sum\limits_{r\ge 1}s_{r}\,=\,\sum\limits_{p\in \mathcal{I}}\sum\limits_{r\ge 1}\sum\limits_{\substack{n\in A_{r} \\ n\equiv 0(\mmod p^{2})}}1\,\le\,
\sum\limits_{p\in \mathcal{I}}\sum\limits_{\substack{n\le x \\ n\equiv 0(\mmod p^{2})}}1\,\le\,\sum\limits_{p>X}\frac{x}{p^{2}}\,<\,xX^{-1},
\]
we find that
\begin{multline*}
\sum\limits_{r\ge 1}|S_{r}|\,<\,\frac{7x}{2}\biggl(\,\frac{1.5}{\sqrt{X\log{X\mathstrut}}}\,+\,\tau(q)\biggl(\frac{\sqrt{q}\log{q}}{x}\biggr)^{\!\! 1/2}\biggr)\log\log{x}
\,+\,\frac{2x}{X}\log\log{x}\,<\\
<\,\frac{7x}{2}\biggl(\,\frac{1.5}{\sqrt{X\log{X}\mathstrut}}\,+\,\tau(q)\biggl(\frac{Y\sqrt{q}\log{q}}{x}\biggr)^{\!\! 1/2}\,+\,\frac{4}{7X}\biggr)\log\log{q}\,<\\
<\,\frac{7x}{2}\biggl(\,\frac{2}{\sqrt{X\log{X}\mathstrut}}\,+\,\tau(q)\biggl(\frac{Y\sqrt{q}\log{q}}{x}\biggr)^{\!\! 1/2}\biggr)\log\log{q}.
\end{multline*}
Finally we get:
\begin{multline*}
|S|\,\le\,(C+3)\,x\,\frac{\log{X}}{\log{Y}}\,+\,\frac{7x}{2}\biggl(\,\frac{2}{\sqrt{X\log{X}\mathstrut}}\,+\,\tau(q)\biggl(\frac{Y\sqrt{q}\log{q}}{x}\biggr)^{\!\! 1/2}\,\biggr)\log\log{q}\,=\\
=\,x\biggl((C+3)\,\frac{\log{X}}{\log{Y}}\,+\,\frac{4\log\log{x}}{\sqrt{X\log{X}\mathstrut}}\,+\,\frac{7}{2}\,\tau(q)\biggl(\frac{Y}{x}\,\sqrt{q}
(\log{q})(\log\log{q})^{2}\biggr)^{\!\! 1/2}\,\biggr).
\end{multline*}
Lemma is proved. $\Box$

Now theorems 1--\,4 easily follow from lemma 5. Indeed, if $x\ge q^{1/2+\vep}$ then we take $X=(\log{q})^{2}$, $Y = q^{\,\vep/4}$ and
suppose $q$ to be so large that $\tau(q)\le q^{\,\vep/9}$. By lemma 5, we get
\begin{multline*}
\Delta_{1}\,\le\,8\vep^{-1}(C+3)\,\frac{\log\log{q}}{\log{q}},\quad \Delta_{2}\,\le\,\frac{7}{\sqrt{2}}\frac{\sqrt{\log\log{q}\mathstrut}}{\log{q}},\\
\Delta_{3}\,\le\,\frac{7}{2}\,q^{\,\vep/9}\biggl(\frac{q^{1/2+\vep/4}(\log{q})}{q^{1/2+\vep}}\biggr)^{\! 1/2}\log\log{q}\,<\,q^{-\vep/4}.
\end{multline*}
Since $8(C+3)<561.7$, then theorem 1 follows.

Further, if $x\ge \sqrt{q}\,e^{(\log\log{q})^{1+\gamma}}$, $\tau(q)\le e^{0.25(\log\log{q})^{1+\gamma}}$ then we set $X=(\log\log{X})^{4+2\gamma}$,
$Y = e^{0.25(\log\log{q})^{1+\gamma}}$. Lemma 5 gives
\begin{multline*}
\Delta_{1}\,\le\,8(2+\gamma)(C+3)\,\frac{\log\log\log{q}}{(\log\log{q})^{1+\gamma}},\quad
\Delta_{2}\,\le\,\frac{7}{\sqrt{2(2+\gamma)}}\,\frac{(\log\log\log{q})^{-0.5}}{(\log\log{q})^{1+\gamma}},\\
\Delta_{3}\,\le\,e^{0.25(\log\log{q})^{1+\gamma}}\biggl(\frac{\sqrt{q}\,(\log{q})e^{0.25(\log\log{q})^{1+\gamma}}}{\sqrt{q}\,e^{(\log\log{q})^{1+\gamma\mathstrut}}}\biggr)^{\!1/2}\,
<\,e^{-0.1(\log\log{q})^{1+\gamma}},
\end{multline*}
and we arrive at theorem 2.

Next, if $x\ge \sqrt{q}\,(\log{q})^{1+2\gamma}$, $\tau(q)\le (\log{q})^{\gamma/4}$ then, setting $X = (\log\log{q})^{4}$, $Y = (\log{q})^{\gamma}$ we obtain
\begin{multline*}
\Delta_{1}\,\le\,4(C+3)\,\frac{\log\log\log{q}}{\gamma\log\log{q}},\quad \Delta_{2}\,\le\,\frac{7(\log\log\log{q})^{-0.5}}{2\log\log{q}},\\
\Delta_{3}\,\le\,(\log{q})^{\gamma/4}\biggl(\frac{\sqrt{q}\,(\log{q})^{1+\gamma}}{\sqrt{q}\,(\log{q})^{1+2\gamma}}\biggr)^{\!1/2}\log\log{q}\,<\,
(\log{q})^{-\gamma/5}.
\end{multline*}
Since $4(C+3)<280.9$ we get the assertion of theorem 3.

Finally, in the case of prime $q$ and $x\ge \sqrt{q}\,(\log{q})e^{(\log\log\log{q})^{1+\gamma}}$ we set $X = (\log\log{q})^{4}$,
$Y = e^{0.5(\log\log\log{q})^{1+\gamma}}$ and then obtain
\begin{multline*}
\Delta_{1}\,\le\,\frac{4(C+3)}{(\log\log\log{q})^{\gamma}},\quad \Delta_{2}\,\le\,\frac{(2\log\log\log{q})^{-0.5}}{\log\log{q}},\\
\Delta_{3}\,\le\,7\biggl(\frac{\sqrt{q}\,(\log{q})e^{0.5(\log\log\log{q})^{1+\gamma}}(\log\log{q})^{2}}{\sqrt{q}\,(\log{q})e^{(\log\log\log{q})^{1+\gamma}}}\biggr)^{\!1/2}
\,<\,e^{-0.1(\log\log\log{q})^{1+\gamma}}.
\end{multline*}
Theorem 4 is proved. $\Box$

We conclude this section with some remarks. Thus, in the case of prime $q$ and $e^{(\log{q})^{4/5}(\log\log{q})^{5}}\le x\le q^{\,4/7}$, the sums $S_{q}(x;f)$ with $f(n) = \mu(n)$, $\mu^{2}(n)$, $\tau_{k}(n)$, $r(n) = \#\{(x,y)\in \mathbb{Z}^{2}\,|\,x^{2}+y^{2}=n\}$ and so on, were estimated in \cite{Korolev_2010}.
Moreover, if $q$ is prime then a slight modification of the arguments from
\cite{Bourgain_Garaev_2014}, \cite{Korolev_2016a} leads to non-trivial bounds
for the very short sums $S_{q}(x;f)$. Namely, the estimate
\[
\sum\limits_{n\le x}f(n)e_{q}(an^{*})\ll xD_{1}^{-1},\quad D_{1}\,=\,\frac{(\log{x})^{3/2}}{\log{q}}\,(\log\log{q})^{-2}
\]
holds for any $x$, $e^{(\log{q})^{2/3}(\log\log{q})^{4/3}}<x\le\sqrt{q}$, and the estimates
\begin{multline*}
\sum\limits_{n\le x}f(n)e_{q}(an^{*})\ll xD_{2}^{-1}\log{D_{2}},\quad \sum\limits_{n\le x}f(n)e_{q}(an^{*}+bn)\ll xD_{2}^{-3/4},\\
D_{2}\,=\,\frac{\log{x}}{(\log{q})^{2/3}}\,(\log\log{q})^{-1/3}
\end{multline*}
holds for $e^{(\log{q})^{2/3}(\log\log{q})^{1/3}}<x\le\sqrt{q}$.

\vspace{0.3cm}

\textbf{\S 5. Some particular cases of the sum $\boldsymbol{S_{q}(x;f)}$}

\vspace{0.5cm}

For some arithmetical functions $f$ and for the prime moduli $q$, the sum $S_{q}(x;f)$
is estimated with power-saving factor for $x\ge q^{1/2+\vep}$. Such estimates
are based on the following very deep result of J.~Bourgain \cite{Bourgain_2005}.

\vspace{0.3cm}

\textsc{Lemma 6.} \emph{Suppose that $q$ is prime, $(a,q)=1$, $0<\vep<0.1$, and let $M, M_{1}$, $N, N_{1}$ satisfy the following conditions:
$q^{\vep}<M,N\le\sqrt{q}$, $M<M_{1}\le 2M$, $N<N_{1}\le 2N$. Suppose also that complex-valued sequences $\alpha_{m}$, $\beta_{n}$ satisfy the inequalities
$|\alpha_{m}|\le \tau_{\ell}(m)$, $|\beta_{n}|\le \tau_{r}(n)$ for some fixed $\ell, r\ge 1$. Then the sum}
\begin{equation}\label{lab_06}
W\,=\,\sum\limits_{M<m\le M_{1}}\sum\limits_{N<n\le N_{1}}\alpha_{m}\beta_{n}e_{q}(am^{*}n^{*}+bmn)
\end{equation}
\emph{obeys the estimate $|W|\le c_{1}MNq^{-c\vep^{4}}$, where $c>0$ is an absolute constant and $c_{1}>0$ depends on $\vep, \ell$ and $r$.}

\vspace{0.3cm}

The original paper of J.~Bourgain \cite{Bourgain_2005} does not contain
the precise expression for the decreasing factor $q^{-c\vep^{4}}$.
It's calculation was made by R.C.~Baker \cite{Baker_2012} for the case $b\equiv 0\pmod{q}$.
However, the arguments of Baker can be adapted without big efforts for the case
$(b,q)=1$.

\vspace{0.3cm}

\textsc{Proof of Theorem 5.} Let us take $y = q^{\,\vep/4}$, $F(n) = e_{q}(an^{*}+bn)$ in the identity
\begin{multline}\label{lab_07}
\sum\limits_{n\le x}\mu(n)F(n)\,=\,-\sum\limits_{k\le y}\mu(k)\sum\limits_{m\le y}\mu(m)\sum\limits_{n\le x/(km)}F(kmn)\,-\\
-\,\sum\limits_{y<m\le x/y}a_{m}\sum\limits_{y<n\le x/m}\mu(m)F(mn)\,+\,2\sum\limits_{n\le y}\mu(n)F(n),
\end{multline}
where
\[
a_{m}\,=\,\sum\limits_{d|m,\;d\le y}\mu(d)
\]
(see, for example, \cite[ch. 2, \S 6, theorem 2]{Karatsuba_Voronin_1992}). Writing the right-hand side of (\ref{lab_07}) in the form
$-\Sigma_{1}-\Sigma_{2}+2\Sigma_{3}$ (where the notations are evident) and using lemma 3, we get:
\[
\bigl|\Sigma_{1}\bigr|\,\le\,\sum\limits_{k,m\le y}4\sqrt{q}\,\log{q}\,\le\,4y^{2}\sqrt{q}\,\log{q}\,<\,xq^{-\vep/4}.
\]
Further, we split the sum $\Sigma_{2}$ to $\ll (\log{q})^{2}$ sums $W$ of the type (\ref{lab_06}), but with the additional restriction
$mn\le x$. Setting $N_{2} = \min{(N_{1}, xm^{-1})}$, $\alpha_{m} = a_{m}$, $\beta_{n} = \mu(n)$ we obtain:
\begin{multline*}
W\,=\,\sum\limits_{M<m\le M_{1}}\sum\limits_{N<n\le N_{1}}\biggl(\frac{1}{q}\sum\limits_{|c|<q/2}\sum\limits_{N<\nu\le N_{2}}e_{q}(c(n-\nu))\biggr)
\alpha_{m}\beta_{n}e_{q}(am^{*}n^{*}+bmn)\,=\\
=\,\sum\limits_{|c|<q/2}(|c|+1)^{-1}W_{c},
\end{multline*}
where the sum $W_{c}$ has the same form as the sum in (\ref{lab_06}), namely:
\begin{multline*}
W_{c}\,=\,\sum\limits_{M<m\le M_{1}}\sum\limits_{N<n\le N_{1}}\alpha(m)\beta(n)e_{q}(am^{*}n^{*}+bmn),\\
 \alpha(m) = \frac{|c|+1}{q}\sum\limits_{N<\nu\le N_{2}}e_{q}(-c\nu),\quad
\beta(n) = e_{q}(cn)\beta_{n}.
\end{multline*}
Obviously, $|\alpha(m)|\le |\alpha_{m}|$, $|\beta(m)|=|\beta_{n}|$. Using lemma 6, we find sequentially
\[
W_{c}\,\ll\,MNq^{-c_{1}\vep^{4}},\quad W\,\ll\,MNq^{-c_{1}\vep^{4}}(\log{q}),\quad \Sigma_{2}\,\ll\,xq^{-c_{1}\vep^{4}}(\log{q})^{3}.
\]
Estimating the sum $\Sigma_{3}$ trivially, we arrive at the desired assertion. $\Box$

\vspace{0.3cm}

Lemma 6 allows one to estimate with power-saving factor and $x\ge q^{1/2+\vep}$ the sums $S_{q}(x;f)$
where the arithmetical function $f$ does not satisfy the condition $|f(n)|\le 1$.

\vspace{0.3cm}

\textsc{Proof of Theorem 6.} In the case $k=1$ the inequality (\ref{lab_08}) follows from lemma 3. Suppose that the estimate (\ref{lab_08}) holds for any sum
$S_{k-1}(y)$ with $q^{1/2+\vep}\le y\le q$ and then verify it for the sum $S_{k}(x)$. Setting $F(n) = e_{q}(an^{*}+bn)$ for brevity and using Dirichlet's hyperbola trick, we find
\begin{multline}\label{lab_08}
S_{k}(x)\,=\,\sum\limits_{uv\le x}\tau_{k-1}(u)F(uv)\,=\,\sum\limits_{u\le \sqrt{x}}\sum\limits_{v\le x/u}\tau_{k-1}(u)F(uv)\,+\\
+\,\sum\limits_{u\le \sqrt{x}}\sum\limits_{v\le x/u}\tau_{k-1}(v)F(uv)\,-\,\sum\limits_{u,v\le\sqrt{x}}\tau_{k-1}(u)F(uv).
\end{multline}
Let $y = q^{\,\vep/2}$. If $u\le y$ ($v\le y$) then $xu^{-1}\ge q^{1/2+\vep/2}$ (correspondingly, $xv^{-1}\ge q^{1/2+\vep/2}$).
By induction, the contribution to $S_{k}(x)$ from the terms with $u\le y$ in the first sum in (\ref{lab_08}), from the terms with $v\le y$ in the second sum and from the terms with $\min{(u,v)}\le y$ in the last sum in (\ref{lab_08}) is estimated  by $O\bigl(xq^{-c\vep^{4}}\bigr)$ in absolute value.

Further, the contribution from the terms with $v\le x^{1/4}$ in the first sum in (\ref{lab_08}) and with $u\le x^{1/4}$ in the second sum does not exceed $x^{3/4}(\log{q})^{k-2}\ll xq^{-1/8}$. Hence,
\[
S_{k}(x)\,=\,\Sigma_{1}+\Sigma_{2}-\Sigma_{3}\,+\,O\bigl(xq^{-c\vep^{4}}\bigr),
\]
where
\begin{multline*}
\Sigma_{1}\,=\,\sum\limits_{y<u\le \sqrt{x}}\;\sum\limits_{x^{1/4}<v\le x/u}\tau_{k-1}(u)F(uv),\quad
\Sigma_{2}\,=\,\sum\limits_{y<u\le \sqrt{x}}\;\sum\limits_{x^{1/4}<v\le x/v}\tau_{k-1}(v)F(uv),\\
\Sigma_{3}\,=\,\sum\limits_{y<u,v\le \sqrt{x}}\tau_{k-1}(u)F(uv).
\end{multline*}
We split every sum $\Sigma_{j}$, $j = 1,2,3$, to $\ll (\log{q})^{2}$ sums $W$ of the type (\ref{lab_06}),
but with the additional restriction $uv\le x$ (in the cases $j = 1,2$). Using the same arguments as above, we get $\Sigma_{j}\ll xq^{-c\vep^{4}}$. Theorem is proved. $\Box$

\renewcommand{\refname}{\normalsize{Bibliography}}

\textsc{Maxim Aleksandrovich Korolev}

Steklov Mathematical Institute of Russian Academy of Sciences

119991, Russia, Moscow, Gubkina str., 8

e\,-mail: \texttt{korolevma@mi.ras.ru}

\end{document}